\documentclass{article}

\usepackage{tikz}

\usepackage{latexsym,amsmath,accents}
\usepackage{amsmath,amsthm,amsfonts,amssymb,graphicx,epsfig,latexsym,color}
\usepackage{epsf,enumerate}
\usepackage{pictex,dcpic}
\usepackage{float}
\usepackage{etex}
\usepackage{wasysym}
\usepackage{ulem}
\usepackage{bbm}
\usepackage{tikz}
\usepackage{tikz-cd}
  \usepackage{tkz-euclide}

\newtheorem{theorem}{Theorem}[section]
\newtheorem{lemma}[theorem]{Lemma}
\newtheorem{corollary}[theorem]{Corollary}
\theoremstyle{definition}

\newtheorem{examples}[theorem]{Examples}

\newtheorem{definition}[theorem]{Definition}

\DeclareMathOperator{\diam}{diam}
\DeclareMathOperator{\asdim}{asdim}
\def\Z{\mathbb Z}
\def\R{\mathbb R}

\title{Groups acting on hyperbolic spaces -- a survey}

\author{Mladen Bestvina}

\begin{document}
\maketitle
\begin{abstract}
This is a (very subjective) survey paper for nonspecialists covering
group actions on Gromov hyperbolic spaces. The first section is about
hyperbolic groups themselves, while the rest of the paper focuses on
mapping class groups and $Out(F_n)$, and the way to understand their
large scale geometry using their actions on various hyperbolic spaces
constructed using projection complexes. This understanding for
$Out(F_n)$ significantly lags behind that of mapping class groups and
the paper ends with a few open questions.
\end{abstract}

\tableofcontents
\section{Introduction}

The goal of this paper is to give a flavor of the developments in
geometric group theory in the last 35 years, focusing on groups acting
on Gromov hyperbolic spaces. The field of geometric group theory is
relatively young and its birth can be attributed to Gromov's paper
\cite{gromov:hypgroups} in 1987, when the subject exploded and
attracted many mathematicians. The term itself was coined by Niblo and
Roller, who organized and named a very influential conference in 1991
\cite{niblo-roller1,niblo-roller2} (though it was possibly used
informally before).  Loosely speaking, geometric group theory studies
groups by looking at their actions on metric spaces and the geometry
and topology of these spaces. Increasingly, methods of other branches
of mathematics, such as dynamics and analysis, are also brought to
bear.

There were, of course, significant developments that can be
comfortably placed within this subject even long before Gromov's
paper. Works of Klein, Dehn, Nielsen, Stallings and others
in some sense form the backbone of the
subject. The theory of groups acting on trees, i.e. Bass-Serre theory
\cite{serre,bass} and its language, will be used freely in these
notes. 
Gromov's celebrated theorem that groups of polynomial growth are
virtually nilpotent \cite{gromov:polygrowth} appeared in 1981, and Gromov's basic philosophy of
viewing groups as metric spaces was eloquently explained in
\cite{gromov:icm}. Of course, the influence on
this subject of the work of Thurston cannot
be overstated. 
Perhaps the development of combinatorial group theory,
focusing on the combinatorics of the words in a finitely presented
group, distracted from a more geometric approach to group theory.

This paper will focus on the part of geometric group theory that
studies groups acting on (Gromov) hyperbolic spaces. In the early days, right
after Gromov's paper, this meant studying (Gromov) hyperbolic
groups. Around 2000, the work of Masur and Minsky \cite{MM1,MM2}
shifted the focus to groups that are not hyperbolic but admit
interesting actions on hyperbolic spaces. The main examples of such
groups are mapping class groups of compact surfaces (the subject of the
papers by Masur and Minsky) and $Out(F_n)$,
the outer automorphism group of a finite rank free group. This survey
will concentrate on these two classes of groups.

The definition of Gromov hyperbolic spaces is modeled on the standard
hyperbolic spaces by ``coarsification'' and captures the fact that
geodesic triangles in the hyperbolic plane are ``thin''. For a
wonderful survey of the history of hyperbolic geometry from
Lobachevsky to 1980 see Milnor's paper \cite{milnor}. For much more
about this subject see Bridson-Haefliger \cite{bridson-haefliger},
Ghys-de la Harpe \cite{Ghys-Harpe} or Dru\c tu-Kapovich
\cite{drutu-kapovich}. 
There are many important topics this survey will not cover,
e.g. relative hyperbolicity \cite{farb}, hyperbolic Dehn filling
\cite{groves-manning, osin}, small cancelation
\cite{arzhantseva, osajda}, uniform embeddings in Hilbert spaces
\cite{sela}, the celebrated work of Agol and Wise, see
e.g. \cite{mb:bulletin}, random walk \cite{maher-tiozzo},
Cannon-Thurston maps \cite{mahan} and many others.

{\bf Acknowledgements.}
I would like to thank all my collaborators over the years, and
particularly Ken Bromberg, Mark Feighn and Koji Fujiwara with whom I
wrote many papers. I had a lot
of fun and learned
many things from you. Here is to the future papers!

\noindent
This work was partially supported by the National Science Foundation,
DMS-1905720.

\section{Hyperbolic groups}

Every finitely generated group $G$ can be viewed as a metric
space. Fix a finite generating set $S$ which is symmetric,
i.e. $S^{-1}=S$. The \textit{word norm} $|g|_S$ of $g\in G$ is the
smallest $n$ such that $g$ can be written as $g=s_1s_2\cdots s_n$ for
$s_i\in S$. Then $d_S(g,h)=|g^{-1}h|_S$ is the \textit{word metric} on
$G$, and
left translations $L_x:g\mapsto xg$ are isometries. More
geometrically, this is the distance function on the vertices of the
Cayley graph $\Gamma_S$, with vertex set $G$, and edges of length 1
between $g$ and $gs$ for $g\in G$ and $s\in S$. If $S'$ is a different
finite symmetric generating set for $G$, the identity map $G\to G$ is
bilipschitz with respect to the two word metrics, and are considered
equivalent. There is a more general equivalence relation between
metric spaces that is very convenient in the subject. Let $(X,d_X)$
and $(Y,d_Y)$ be metric spaces. A (not necessarily continuous)
function $f:X\to Y$ is a \textit{quasi-isometry} if there is a number
$A>0$ such that
$$\frac 1A d_X(a,b)-A\leq d_Y(f(a),f(b))\leq A~ d_X(a,b)+A$$ for all
$a,b\in X$, and every metric ball of radius $A$ in $Y$ intersects the
image of $f$. Without the second condition, $f$ is a
\textit{quasi-isometric embedding} (when we want to refer to the
constant $A$ we say $A$-quasi-isometric embedding). Two metric spaces
are \textit{quasi-isometric} if there is a quasi-isometry between
them, and this is an equivalence relation. For example, inclusion
$\Z\hookrightarrow\R$ is a quasi-isometry, as is any bilipschitz
homeomorphism or a finite index inclusion between finitely generated
groups equipped with word metrics. More generally, the following is
considered to be the Fundamental Theorem of Geometric Group Theory.

\begin{theorem} (Milnor \cite{milnor-qi}, \v{S}varc \cite{svarc})
  Suppose a group $G$ acts properly and cocompactly by isometries on a
  proper geodesic metric space $X$. Then $G$ is finitely generated and
  any orbit map $G\to X$ is a quasi-isometry.
\end{theorem}

A metric space is proper if closed metric balls are compact, and it is
\textit{geodesic} if any two distinct points $a,b$ are joined by a subset
isometric to the closed interval $[0,d(a,b)]$. For example, cocompact
lattices in a simple Lie group are quasi-isometric to each other. The
``Gromov program'' is to classify groups, at least in a given class,
up to quasi-isometry.

According to Gromov, the following definition was given by Rips. There
are several other definitions, all of which are equivalent up to
changing the value of $\delta$, see
\cite{bridson-haefliger,drutu-kapovich}.

\begin{definition}
  Let $\delta\geq 0$. A geodesic metric space $X$ is
  \textit{$\delta$-hyperbolic} if in any geodesic triangle each side is
  contained in the $\delta$-neighborhood of the other two sides. We
  say $X$ is \textit{hyperbolic} if it is $\delta$-hyperbolic for some
  $\delta\geq 0$.
\end{definition}

\begin{figure}[h]
  \begin{center}
\begin{tikzpicture}[scale=0.6,y=-1cm]

\fill[fill={rgb:black,2;white,10}] plot [smooth cycle] coordinates
     {(13.9,15.1)  (14.9,15.4) (15.2,14.4) (14.3,13.2) (13.1,11.5)
       (12.2,9.6) (12.1,7.6) (11.5,6.5) (10.4,7.1) (10.1,8.4)
       (10,10.1) (10.3,11.7) (11.5,13.4) (12.8,14.5)};
\fill[fill={rgb:black,2;white,10}] plot [smooth cycle] coordinates {(11.3,6.5)
  (10.6,6.9) (10.1,8.2) (9.6,10) (9.1,10.9) (8.2,11.5) (6.9,11.8)
  (5.7,12.2) (5.7,13.1) (6.2,13.5) (7.2,13.3) (8.2,13.2) (9.5,13)
  (10.5,12.6) (11.3,12) (11.8,10.9) (11.9,9.5) (12,7.7) (11.6,6.7)};
\draw[black] plot [smooth] coordinates {(6.2,12.7) (10,11.3) (11.2,7)};
\draw[black] plot [smooth] coordinates {(6.2,12.7) (10.3,12.22) (14.3,14.5)};
\draw[black] plot [smooth] coordinates {(11.2,7) (11.2,11)
  (14.3,14.5)};
\node[xshift=-0.3cm] at (6.2,12.7) {A};
\node[xshift=0.3cm] at (14.3,14.5) {B};
\node[yshift=0.3cm] at (11.2,7) {C};

\end{tikzpicture}%
\caption{The union of the $\delta$-neighborhoods of two sides contains
  the third.}
\end{center}
  \end{figure}
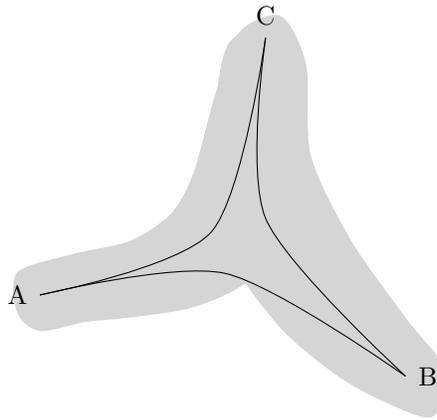

For example, trees are 0-hyperbolic and so are complete
simply-connected Riemannian manifolds of sectional curvature
$\leq-\epsilon<0$. A fundamental property of hyperbolic spaces is the
Morse Lemma, proved by Morse \cite{morse}, Busemann \cite{busemann}
and Gromov \cite{gromov:hypgroups} in increasing generality.

\begin{lemma}[Morse Lemma]
  There is a number $D=D(\delta,A)$ such that for any
  $\delta$-hyperbolic space $X$ and any $A$-quasi-isometric
  embedding $f:[a,b]\to X$ the image of $f$ is contained in the
  $D$-neighborhood of any geodesic from $f(a)$ to $f(b)$.
\end{lemma}

It then quickly follows that if two geodesic spaces are
quasi-isometric and one is hyperbolic, so is the other. In particular,
groups that act properly and cocompactly by isometries on proper
hyperbolic spaces are hyperbolic. 

Hyperbolic groups are well-behaved, both topologically and
geometrically, and they are generic, so they form a model class of
groups in geometric group theory. We now elaborate.

\subsection{Classification of elements}
Let $G$ be a hyperbolic group. If $g\in G$ has finite order, then
there is a coset $\langle g\rangle x$ that has diameter $\leq
4\delta+2$, so in particular there is an a priori bound on the order
in terms of $\delta$ and the number of generators. This is proved by a
coarse version of the standard argument that a bounded set in $\R^n$
(or any Hadamard manifold) is contained in a unique closed ball of
smallest radius. If $g$ has infinite order, then $k\mapsto g^kx$ is a
quasi-isometric embedding for every $x\in G$, and $g$ is 
  \textit{loxodromic}. 

\subsection{The Rips complex}
The classical Cartan-Hadamard theorem states that closed manifolds of
nonpositive sectional curvature have contractible universal cover. In
a similar way, every hyperbolic group $G$ acts properly and
cocompactly on a contractible simplicial complex, called the \textit{Rips
  complex}. It is constructed as follows. Fix a number $d>0$ and form
the complex $P_d(G)$: the set of vertices is $G$, and a set
$\{v_0,v_1,\cdots,v_n\}$ of distinct vertices forms a simplex if
$d(v_i,v_j)\leq d$ for all $i,j$. This is a version of the Vietoris
approximation of a metric space by a simplicial complex, except here
we think of $d$ as being large.

\begin{theorem}
  For $d>4\delta+6$ $P_d(G)$ is contractible.
\end{theorem}

So for example if $G$ is torsion-free, the quotient $P_d(G)/G$ is a
finite classifying space for $G$, and in any case $G$ is finitely presented,
and has a classifying space with finitely many cells in each
dimension. Every finite subgroup of $G$ fixes a point of $P_d(G)$ (for
$d$ large), so it follows that $G$ has finitely many conjugacy classes
of finite subgroups. Interestingly, it is not known whether every
infinite hyperbolic group is virtually torsion-free, or even if it
always has a proper subgroup of finite index.

\subsection{Subgroups}
If $g\in G$ has infinite order, there is a unique maximal virtually
cyclic subgroup $E(g)$ of $G$ that contains $g$, and $E(g)$ also
contains the normalizer of $g$. It follows that $G$ cannot contain
$\Z^2$ as a subgroup. Translation length considerations show that $G$
cannot contain any Baumslag-Solitar groups $B(m,n)=\{a,t\mid
ta^mt^{-1}=a^n\}$, $m,n\neq 0$, as subgroups. The long standing open
question whether every group with finite classifying space and not
containing any $B(m,n)$ is necessarily hyperbolic was recently
answered in the negative \cite{bruno}.

\subsection{Boundary}
Inspired by the visual boundary of Hadamard manifolds, Gromov defined
a boundary $\partial G$ of a hyperbolic group (or a proper geodesic
metric space which is hyperbolic). It is a compact
metrizable space and a point is represented by a quasi-geodesic ray
$\Z_+\to G$, with two rays representing the same boundary point if
their images stay a bounded distance apart. The topology is based on
the principle that rays issuing from a basepoint and with fixed
quasi-geodesic constants will stay longer together if they represent
points that are closer together. If $G$ is infinite and virtually cyclic then
$\partial G$ consists of two points, and if $G$ is not
virtually cyclic (termed ``non-elementary'') $\partial G$ has no
isolated points.

There is also a natural topology on
the union $$\overline X=P_d(G)\sqcup\partial G$$ of the Rips complex and the
Gromov boundary
that makes it into a compact metrizable space, and $G$ acts naturally by
homeomorphisms. Loxodromic elements act by north-south dynamics on
$\overline X$. The most important property of the boundary, used for
example in the proof of Mostow rigidity \cite{mostow}, is the
following:

\begin{theorem}
  Let $f:X\to Y$ be a quasi-isometry between two hyperbolic proper geodesic
  mmetric spaces. Then $f$ extends to a homeomorphism $\partial X\to\partial
  Y$.
\end{theorem}

\begin{theorem}\cite{b-mess}
  $\overline X$ is a Euclidean retract, i.e. it is contractible,
  locally contractible and finite-dimensional. The covering dimension
  of $\partial G$ can be computed from the cohomology of $G$ and in
  particular, if $G$ is torsion-free, $\dim \partial G$ equals the
  cohomological dimension of $G$ minus 1, and in any case the rational
  cohomological dimension of $\partial G$ equals the rational cohomological
  dimension of $G$ minus 1.
\end{theorem}

\subsection{Asymptotic dimension}
In \cite{gromov:asymptotic} Gromov introduced many quasi-isometry
invariants of groups and spaces. Here we focus on \textit{asymptotic
  dimension}. Let $X$ be any metric space. For an integer $n\geq 0$ we
write $asdim(X)\leq n$ provided that for every $R>0$ there exists a
cover of $X$ by uniformly bounded sets such that every ball of radius
$R$ in $X$ intersects at most $n+1$ elements of the cover. This is the
``large scale'' analog of the usual covering dimension. For example,
$asdim(\R^n)=n$ and $asdim(T)\leq 1$ for a tree $T$ with the geodesic
metric. This is a quasi-isometry invariant, so it is well defined for
finitely generated groups as well. See \cite{bell-dranishnikov} for
the basic properties of $asdim$. There are many groups that contain
$\Z^n$ for every $n$, and they will have infinite asymptotic
dimension. However, Gromov proved:

\begin{theorem}\cite{gromov:asymptotic}
  Every hyperbolic group has finite asymptotic dimension.
\end{theorem}

One can hardly make a claim that one understands the large scale
geometry of a group if its asymptotic dimension is not known to be
finite or infinite.  However, the significance of the theorem became
particularly clear with the work of Guoliang Yu \cite{yu} (see also
\cite{dranishnikov-ferry-weinberger}), who proved that groups with
finite $asdim$ and finite classifying space satisfy the Novikov
conjecture (this predicts the possible placement of Pontrjagin classes
in the cohomology ring of a closed oriented manifold with the given
fundamental group).

An even stronger conjecture in manifold topology is the Farrell-Jones
conjecture. If it holds for a (torsion-free) group $G$ then one can in
principle compute the set of closed manifolds homotopy equivalent to a
given closed manifold of dimension $\geq 5$ and fundamental group
$G$. Following the work of Farrell and Jones, there has been a great
progress in proving the Farrell-Jones conjecture for many groups. For
hyperbolic groups this was done by Bartel, L\"uck and Reich
\cite{bartels-lueck-reich}, see also \cite{bartels} for a proof using
coarse methods that generalize to other groups.

\subsection{JSJ decomposition}
For simplicity, we now assume that $G$ is a torsion-free hyperbolic
group. By Grushko's theorem \cite{grushko,stallings:grushko} $G$ can
be decomposed as a free product $G=G_1*G_2*\cdots *G_k*F_r$ where each
$G_i$ is noncyclic and freely indecomposable and $F_r$ is a free
group. Each $G_i$ is a 1-ended group by the celebrated theorem of
Stallings \cite{stallings}, meaning that the Cayley graph of $G_i$ has
one end (every finite subgraph has only one unbounded complementary
component). Quite unexpectedly, Rips-Sela \cite{rips-sela:jsj} discovered
a further structure theorem for 1-ended hyperbolic groups (the theorem
applies to many groups that are not hyperbolic as well). The theorem
is motivated by the Jaco-Shalen-Johanssen torus decomposition theorem
for 3-manifolds, which provides a canonical decomposition of an
aspherical closed orientable 3-manifold by cutting along pairwise
disjoint tori so that each piece either has many tori (it is Seifert
fibered), or it is not an $I$-bundle and has no essential tori (except
on the boundary, and then by Thurston's hyperbolization theorem it is
hyperbolic), or it is an $I$-bundle. The Rips-Sela theorem can be
stated as follows:

\begin{theorem}
  Let $G$ be a 1-ended torsion-free hyperbolic group. Then $G$ is a
  finite graph of groups with all edge group infinite cyclic, and with
  vertex groups $V$ coming in three types:
  \begin{enumerate}
    \item [(QH)] $V$ is the fundamental group of a compact surface (with
      a pair of intersecting 2-sided simple closed curves) and
      the incident edge groups correspond exactly to the boundary
      components, 
      \item [(rigid)] $V$ is not cyclic and does not admit a nontrivial
        splitting over a cyclic group such that all incident edge
        groups are elliptic, and
        \item [(cyclic)] $V$ is cyclic.
  \end{enumerate}
\end{theorem}

See also \cite{dunwoody-sageev, papa-fu, gl:jsj} for different proofs
and generalizations, and \cite{bowditch:jsj} for how to read off the
JSJ decomposition purely from the boundary of $G$. For example, a
splitting over $\Z$ gives a pair of points in $\partial G$ that
together separate $\partial G$, and Bowditch shows how to go in the
other direction.
Thus the QH
vertices give rise to many splittings of $G$ over cyclic groups (one
for every simple closed curve), while
rigid vertices give rise to none.

\begin{figure}[h]
  \begin{center}
\begin{tikzpicture}[scale=0.4,y=-1cm]
\draw[dashed,black] (13.13776,8.82143) +(147:0.87938) arc (147:41:0.87938);
\draw[black] (5.55,11.85) +(-114:2.35053) arc (-114:-66:2.35053);
\draw[black] (13.03857,10.21) +(-125:1.1117) arc (-125:-47:1.1117);
\draw[black] (9.15,7.55) +(-58:0.65192) arc (-58:58:0.65192);
\draw[dashed,black] (9.85,7.55) +(-122:0.65192) arc (-122:-238:0.65192);
\draw[fill={rgb:black,2;white,10}] (4.2,6.4) rectangle (7.1,8.8);
\draw[fill={rgb:black,2;white,10}] (11.8,6.4) rectangle (14.7,8.8);
\draw[black] (7.1,7) -- (11.8,7);
\draw[black] (7.1,8.1) -- (11.8,8.1);
\draw[black] plot [smooth] coordinates {(6,8.8) (9.2,13) (12.7,8.8)};
\draw[black] plot [smooth] coordinates {(4.7,8.8) (5.2,16.1) (12.7,16) (13.8,8.8)};
\draw[black] plot [smooth] coordinates{(7.8,14.8) (8.5,15.8) (9.2,14.9)};
\draw[black] plot [smooth] coordinates {(8,15.2) (8.5,14.8) (9,15.2)};
\draw[dashed,black] plot [smooth] coordinates{(4.6,9.7) (5.5,10) (6.5,9.7)};

\end{tikzpicture}%
\caption{A possible JSJ decomposition of a group $G$, with two rigid
  vertices and one QH vertex.}
  \end{center}
  \end{figure}
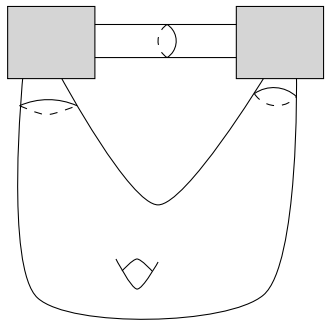

We can picture $G$ as the
fundamental group of the space obtained from a disjoint union of
compact surfaces, ``black boxes'' and circles by attaching cylinders
according to the graph of groups. The JSJ decomposition is not quite
unique, but there are standard moves that transform one such
decomposition to another. For example, sometimes one can slide one
cylinder over another if they meet at a common circle. The main
feature of a JSJ decomposition is that splittings over cyclic groups
can be ``read off'', at least up the standard moves, just like all
essential tori in a 3-manifold can be read off from its JSJ
decomposition. 

\subsection{The Combination theorem}
This is also motivated by 3-manifold theory. The classical
Klein-Maskit combination theorem gives conditions under which two
discrete groups $A,B$ of isometries of hyperbolic space $\mathbb H^3$
with intersection $C=A\cap B$ generate the amalgam $A*_C
B$. Thurston's Hyperbolization Theorem \cite{thurston,morgan} is
proved by cutting the 3-manifold into pieces, and then inductively
constructing a hyperbolic structure when gluing the pieces
together. There are two opposite extremes in the kinds of gluings,
when the intersection of the pieces is quasi-isometrically embedded on
both sides, and when it is exponentially distorted. The latter arises
when the 3-manifold fibers over the circle and the monodromy is
pseudo-Anosov. The following is the hyperbolic group analog.

\begin{theorem}[\cite{bf:combination,bf:combination2}]
  Let $G$ be the fundamental group of a finite graph of hyperbolic
  groups so that each edge group is quasi-isometrically embedded in
  both vertex groups (but not necessarily in $G$). Assume the ``annuli
  flare'' condition. Then $G$ is a hyperbolic group.
\end{theorem}

The precise definition of the annuli flare condition is a bit
technical, but let us mention two special cases. The first is when the
graph of groups is \textit{acylindrical}, that is, for some $M>0$ the
stabilizer of every segment of length $M$ in the associated Bass-Serre
tree is finite. In this case there are no (long) annuli at all. The
other case is that of a hyperbolic automorphism $\phi:H\to H$ of a
hyperbolic group $H$. This means that there is $M>0$ such that for
every element $h\in H$ of sufficiently large word length $|h|$ we have
$$\max\{|\phi^M(h)|,|\phi^{-M}(h)|\}\geq 2|h|$$
so in this case the induced infinite annulus defined on $S^1\times\R$
sending $S^1\times\{K\}$ to the loop determined by $\phi^K(h)$ flares
exponentially. 
Aside from automorphisms of closed surface groups induced by
pseudo-Anosov homeomorphisms, there are many examples (in fact, they
are generic in the sense of random walk
\cite{kapovich-maher-pfaff-taylor}) of hyperbolic
automorphisms of free groups coming from train track theory
\cite{b-handel}. The combination theorem then implies that the mapping
torus $H\rtimes_\phi \Z$ is hyperbolic.

The combination theorem has also been used to study hyperbolicity of
extensions of free or surface groups in terms of the monodromy
homomorphism from the quotient group to the mapping class group or
$Out(F_n)$, giving rise to \textit{convex cocompact subgroups} of
these groups \cite{farb-mosher, ursula, kent-leininger,
  dowdall-taylor:cc}. 

\subsection{Random groups are hyperbolic}
The most straightforward way to talk about ``random groups'' is the
following model. Fix integers $k\geq 2$ and $m\geq 1$ and for integers
$n_1,\cdots,n_m$ consider the finite set
$$N(k,m;n_1,\cdots,n_m)$$
of all group presentations with $k$ generators and $m$ relators of
lengths $n_1,\cdots,n_m$. Say that a \textit{random group has property P} if
the fraction of groups in $N(k,m;n_1,\cdots,n_m)$ that have P goes to
1 as $\min\{n_1,\cdots,n_m\}\to\infty$.

\begin{theorem}\cite{olshanski,champetier}
  A random group is hyperbolic and its boundary is the Menger curve.
\end{theorem}

Thus a random group has rational cohomological dimension 2 and does not split
over a finite or a 2-ended group.

Gromov \cite{gromov:random} introduced a more sophisticated random
model for groups, called the \textit{density model}, that depends on a
parameter $d\in (0,1)$ and properties of random group depend on the
chosen range of $d$. For more information see \cite{ghys:random,
  ollivier}. 

\subsection{$\R$-trees and applications}
$\R$-trees are metric spaces such that any two distinct points $x,y$ are
contained in a unique subspace homeomorphic to a closed interval in $\R$
with $x,y$ corresponding to the endpoints, and this subspace is
isometric to a closed interval. Simplicial trees with the
length metric induced by identifying edges with closed intervals are
examples of $\R$-trees. More generally, $\R$-trees can have a dense
set of ``vertices'' (points whose complement has more than two
components). For example, let $T=\R^2$ as the underlying set, and
define the metric $d$ as follows: $d(x,y)=|x-y|$ is the Euclidean
distance if $x,y$ are on the same vertical line, and otherwise if
$x=(x_1,x_2)$, $y=(y_1,y_2)$, then
$d(x,y)=|y_1|+|y_2|+|x_1-x_2|$. Thus one imagines train lines running
along all vertical lines and along the $x$-axis, with the distance
function being the shortest train trip.

$\R$-trees were put to good use by Morgan and Shalen
\cite{ms1,ms2,ms3} in their work on hyperbolization of 3-manifolds
following Thurston's work.

The importance of $\R$-trees in geometric group theory comes from two
principles that we briefly review. Let $X$ be a proper hyperbolic
space with the isometry group of $X$ acting with coarsely dense
orbits. 
\begin{enumerate}
  \item A sequence of actions of a finitely generated group $G$ on
    $X$ either, after a subsequence, converges
    (after conjugations) to an isometric action on $X$, or else it
    converges to an isometric action on an $\R$-tree.
    \item There is a theory analogous to the Bass-Serre theory, called
      the ``Rips machine'', that explains the structure of a group
      acting isometrically on an $\R$-tree
      from the stabilizers of the action (under some technical
      conditions).
\end{enumerate}

\subsection{Hyperbolic spaces degenerate to $\R$-trees}\label{2.10}
This construction is due to F. Paulin \cite{paulin} and the author
\cite{mb:degenerations}. See also \cite{handbook}.
We fix a group $G$ and a finite generating set
$a_1,\cdots,a_n$. Suppose we are given
an isometric action $\rho:g\mapsto \rho(g):X\to X$ of $G$ on a
proper $\delta$-hyperbolic space $X$, defined up to conjugation by an
isometry of $X$. We impose the mild assumption that the action is
\textit{nonelementary}, i.e. the function
$$x\mapsto \max_j\{d_X(x,a_j(x))\}$$ is a proper function $X\to
       [0,\infty)$. We then choose a basepoint $x_\rho\in X$ where the
         minimum is attained. Identifying $G$ with the orbit of
         $x_\rho$, this induces a left-invariant (pseudo)metric on $G$
         via
         $$d_\rho(g,h)=d_X(g(x_\rho),h(x_\rho))$$

This metric is ``hyperbolic'', although $G$ as a discrete set is not a
geodesic metric space.  To make this precise, it is convenient to give
Gromov's reformulation of $\delta$-hyperbolicity, in terms of the
``4-point condition''. For $a,b\in X$ define the ``Gromov
product'' $$(a\cdot b)=\frac 12
(d_X(x_\rho,a)+d_X(x_\rho,b)-d_X(a,b))$$ Thus when $X$ is a tree
$(a\cdot b)$ is the distance between $x_\rho$ and $[a,b]$, and in
general it is within $2\delta$ of it. If $a,b,c\in X$ then consider
the 3 numbers $(a\cdot b),(b\cdot c), (c\cdot a)$. When $X$ is a tree,
the two smaller numbers are equal. Gromov's 4-point condition is that
the two smaller numbers are within $\delta$ of each other. Up to
changing the value of $\delta$, a geodesic metric space is hyperbolic
if and only if it satisfies the 4-point condition. Moreover, if the
4-point condition holds with $\delta=0$, then the space can be
isometrically embedded in an $\R$-tree.

Returning to our setup, assume now that $\rho_i$ is a sequence of
isometric actions of $G$ on $X$, $x_{\rho_i}$ are the corresponding
basepoints, and $d_{\rho_i}$ the induced metrics on $G$. They all
satisfy the 4-point condition with a fixed $\delta$. There are now two
cases, up to passing to a subsequence. Define
$D_i=\max_j\{x_{\rho_i},a_j(x_{\rho_i})\}$. 

\textbf{Case 1.} $D_i\to\infty$. Then rescale the metrics $d_{\rho_i}$ by
$D_i$, i.e. consider $d_{\rho_i}/D_i$. After a subsequence, this will
converge to a (pseudo)metric on $G$ which will now satisfy the 4-point
condition with $\delta=0$. Thus $(G,d)$ can be isometrically embedded
into a (unique) $\R$-tree $T$ and there will be an induced isometric
action of $G$ on $T$. Thanks to the careful choice of basepoints, this
action will not have a global fixed point. 

\textbf{Case 2.} $D_i$ stays bounded. Under the mild condition that the
isometry group of $X$ acts with coarsely dense orbits, we can
conjugate the given actions so that all $x_{\rho_i}$ belong to a fixed
bounded set. Since $X$ is proper, there is a further subsequence so
that $\rho_i$ converge to an isometric action $\rho$ of $G$ on $X$. 

\subsection{The Rips machine}
If a group acts freely on a simplicial tree, it is necessarily
free. This simple instance of Bass-Serre theory follows quickly from
covering space theory. However, this is not true for $\R$-trees. For
example, $\Z^n$ acts freely on $\R$ by letting basis elements act by
$n$ rationally independent translations. More interestingly, closed
surfaces of Euler characteristic $<-1$ admit measured foliations with
simple singularities and with all leaves being trees (and all but
finitely many are lines), see \cite{thurston:bulletin}. Lifting to the
universal cover, the transverse measure turns the leaf space to an
$\R$-tree and the deck group induces a free action of the fundamental
group of the surface on this $\R$-tree.

Suppose now we are given an isometric action of a finitely presented
group $G$ on an $\R$-tree $T$. We make a technical condition that the
action is \textit{stable} meaning that for every arc $I\subset T$
there is a subarc $J\subset I$ such that the stabilizer of $J$ is
equal to the stabilizer of any further subarc of $J$. This property is
frequently satisfied for actions on $\R$-trees obtained by
degenerating $\delta$-hyperbolic spaces described above. We then fix a
finite simplicial 2-complex $K$ with $G=\pi_1(K)$ and construct a
$G$-equivariant map $\tilde K\to T$, called a \textit{resolution} of
$T$. Point inverses form a foliation of $\tilde K$ (with certain
standard singularities) which descends to $K$. The Rips machine
transforms $K$ with this foliation, without changing the fundamental
group nor the fact that the universal cover resolves $T$, and puts it
in a certain ``normal form''. The pieces of this normal form are
foliated subcomplexes that occur, very surprisingly, in only the
following four types:
\begin{enumerate}
  \item [(simplicial)] leaves are compact and the piece resolves a
    simplicial tree,
    \item [(surface)] the piece is a surface (perhaps with boundary)
      and the nonboundary leaves are trees as above,
    \item [(axial)] the piece resolves the tree which is a line, and
    \item [(Levitt)] the piece is of \textit{Levitt type}.
\end{enumerate}

Levitt type foliations were first constructed by G. Levitt
\cite{levitt:sierpinski}. Generic leaves are 1-ended graphs, and in
fact they are quasi-isometric to 1-ended trees with finite graphs
attached. In addition to proving this classification, the Rips machine
also provides the structure of the group corresponding to these cases,
and particularly in the Levitt case. It turns out that if there is a
Levitt piece then $G$ always splits along a subgroup which fixes an
arc in $T$. The other three cases are classical, with the simplicial
case amounting to Bass-Serre theory. As an example, Rips proved the
conjecture of Morgan and Shalen that any finitely generated group
acting freely on an $\R$-tree is isomorphic to the free product of
surface groups and free abelian groups.
For more details see \cite{folpap,french}.

\subsection{Applications}
We mention some of the applications of $\R$-trees; for more see
\cite{handbook}. They are a basic tool in the theory of
$Out(F_n)$. Zlil Sela used them extensively in his seminal work on the
Tarski problems \cite{sela1,sela2,sela4,sela51,sela52,sela7,sela8}.

\subsubsection{Automorphisms of hyperbolic groups.}
Let $G$ be a 1-ended hyperbolic group and for simplicity we assume it
is torsion-free. Combining Paulin's construction
\cite{paulin:out} with the Rips machine we get:

\begin{theorem}
  If $G$ does not split over $\Z$ then $Out(G)$ is finite.
\end{theorem}

This is analogous to a consequence of Mostow Rigidity that $Out(G)$ is
finite when $G$ is the fundamental group of a closed hyperbolic
$n$-manifold with $n\geq 3$.

The proof goes like this. Assuming $Out(G)$ is infinite, choose a
sequence $f_i$ of automorphisms in distinct classes and consider
isometric actions $\rho_i$ of $G$ on itself given by left translations
twisted by $f_i$, i.e. $g\mapsto (h\mapsto f_i(g)h)$. Since $f_i$ are
distinct in $Out(G)$ we see that we are in Case 2 of the construction
outlined above and we obtain an isometric action of $G$ on an
$\R$-tree and with arc stabilizers cyclic (or trivial). The Rips
machine now yields a splitting of $G$ over a cyclic group.

A proper generalization of this theorem was given by Z. Sela. Fix a
JSJ decomposition of $G$. There are now ``visible'' automorphisms of
$G$ realized as compositions of powers of Dehn twists in the cylinders
and homeomorphisms of the QH vertices, which
are surfaces.

\begin{theorem}\cite{rips-sela:sr1}
  The subgroup of visible automorphisms has finite index in $Out(G)$.
\end{theorem}

The proof is quite a bit harder. The idea is that if the index is
infinite one can choose a sequence of automorphisms $f_i$ in distinct
cosets of the visible subgroup. In addition, one chooses the $f_i$'s
to be the ``shortest'' in their cosets. Then one argues that the
action in the limit produces a ``new'' splitting of $G$, one not
explained by the JSJ, or else the $f_i$ could be shortened for large
$i$.

Recall that a group $G$ is \textit{Hopfian} if every surjective endomorphism of
$G$ is an automorphism and it is co-Hopfian if every injective
endomorphism is an automorphism. For example, nontrivial free groups
are not co-Hopfian. By adapting the above methods to endomorphisms,
Sela proved:

\begin{theorem}\cite{sela:hopf,sela:sr2}
  Let $G$ be torsion-free hyperbolic. Then $G$ is Hopfian. If $G$ is
  1-ended it is also co-Hopfian.
\end{theorem}

In 1911 Max Dehn proposed three algorithmic problems about groups: the
word problem (decide if a word in the generators represents the
trivial element), the conjugacy problem (decide if two words in the
generators represent conjugate elements) and the isomorphism problem
(decide if two groups given by presentations are isomorphic). Dehn
solved the word problem for surface groups and his solution
generalizes to hyperbolic groups. There is also a similar solution of
the conjugacy problem for hyperbolic groups, see
\cite{gromov:hypgroups}. The isomorphism problem takes more work and
uses $\R$-trees. For torsion-free hyperbolic groups that don't split
over cyclic subgroups the isomorphism problem was solved by Sela
\cite{sela:isomorphism}, and for general hyperbolic groups by
Dahmani-Guirardel \cite{dahmani-guirardel}. 

Even though hyperbolic groups are generally very well behaved, they
also contain a certain amount of pathologies, see
e.g. \cite{bridson}. 

\subsubsection{Local conectivity of $\partial G$}
The use of $\R$-trees completed the proof of the following theorem.

\begin{theorem}
  If $G$ is a 1-ended hyperbolic group, then $\partial G$ is locally
  connected (as well as connected).
\end{theorem}

There are several ingredients in the proof. First, \cite{b-mess} shows
that if $\partial G$ is not locally connected then it has (many) cut
points. Bowditch \cite{bowditch:jsj} then shows that $G$ acts on an
$\R$-tree constructed as a kind of a ``dual'' tree, which doesn't come
with a metric but can be endowed with one using
\cite{levitt:nonnesting}. The Rips machine then yields a splitting of
$G$ over a 2-ended group, finishing the proof if such splittings do
not exist. Swarup \cite{swarup} finished the proof in the general case
by showing how to continue refining these splittings (in the presence
of cut points in $\partial G$) until the full
JSJ decomposition is obtained, at which point a contradiction arises
with any further splitting.

\subsubsection{Thurston's compactness theorem.}
With the machinery of $\R$-trees one can give a quick proof of the
following theorem.

\begin{theorem}\cite{thurston:compactness}
  Let $M$ be a compact aspherical 3-manifold whose fundamental group does not
  split over a cyclic group. Then the space of hyperbolic structures $H(M)$
  on $M$ is compact.
\end{theorem}

The space $H(M)$ is the space of discrete and faithful representations
of $G=\pi_1(M)$ into the orientation isometry group $PSL_2(\mathcal
C)$ of hyperbolic 3-space $\mathbb H^3$, up to conjugacy (it takes
some work to see that the quotient of $\mathbb H^3$ by such a group is
homeomorphic to the interior of $M$). Indeed, to rule out Case 2
above, one shows that the limiting action on an $\R$-tree is stable
and has abelian arc stabilizers (which follows from discreteness and
faithfulness).

\section{Mapping class groups}

A fundamental shift in the subject occurred after the work of Masur
and Minsky \cite{MM1,MM2} on mapping class groups, the work that set
the foundations for an eventual understanding of the large scale
geometry of these groups. Mapping class groups are not hyperbolic
(except for some sporadic surfaces) but naturally act on hyperbolic
spaces.

We start by recalling some definitions. Let $S$ be an orientable
surface of finite type, i.e. one obtained from a closed orientable
surface by removing finitely many points (called punctures). The group
$Homeo_+(S)$ of orientation preserving homeomorphisms of $S$ has the
natural compact-open topology which makes it locally path-connected,
and the mapping class group (or the Teichm\"uller modular group)
$Mod(S)$ is the discrete group of (path) components of
$Homeo_+(S)$. Classically, this group has been studied since the early
20th century. A very nice introduction to the subject is the book
\cite{farb-margalit}, and we will freely use the standard
concepts. For example, the subgroup $PMod(S)$ of ``pure'' mapping
classes (those that fix the punctures) is generated by finitely many
Dehn twists and the group will not be hyperbolic if $S$ is big enough
to contain two essential (not bounding a disk or a punctured disk)
nonparallel (not cobounding an annulus) disjoint simple closed curves.

To the surface $S$ Harvey \cite{harvey} associates a simplicial
complex $\mathcal C=\mathcal C(S)$, called the \textit{curve complex}
of $S$. A vertex is an isotopy class of essential simple closed
curves. A collection of distinct vertices spans a simplex if each pair
can be represented by curves that intersect minimally (most of the
time this means ``disjointly'', but in a torus punctured at most once
it means ``once'' and in a four times punctured sphere it means
``twice''). For the purposes of this discussion we restrict to the
1-skeleton (called the \textit{curve graph}), which we equip with the
length metric with all edges of length 1. The group $Mod(S)$ acts
naturally on $\mathcal C(S)$. For some very small surfaces, like a 3
times punctured sphere, the curve complex is empty, but otherwise it
is infinite, and even locally infinite, a big contrast with Cayley
graphs of hyperbolic groups. In a similar way, one can define the arc
complex of a surface with punctures.

\begin{theorem}\cite{MM1}
  $\mathcal C(S)$ is hyperbolic. An element of $Mod(S)$ acts
  loxodromically if and only if it is pseudo-Anosov. 
\end{theorem}

Here are some ideas in the original proof, which uses Teichm\"uller
theory. Let $\mathcal T=\mathcal T(S)$ be the Teichm\"uller space of
$S$, i.e. the space of all (marked) hyperbolic structures on
$S$. There is a natural coarse map $\pi:\mathcal T\to\mathcal C$ that
to a hyperbolic metric on $S$ assigns (the isotopy class of) a
shortest simple closed geodesic. Any two points in $\mathcal T$ are
joined by a unique Teichm\"uller geodesic, and their images under
$\pi$ form a family of coarse paths in $\mathcal C$ satisfying (and
this needs proof):
\begin{itemize}
  \item any two points in $\mathcal C$ are connected by some
    such path,
    \item the family is closed under taking subpaths,
    \item any two paths in the collection starting at nearby points
      are contained in each other's uniform Hausdorff neighborhood
      (i.e. they \textit{fellow travel}),
      and
    \item triangles formed by these paths are uniformly thin.
\end{itemize}

Thus the collection behaves like the collection of geodesics in a
hyperbolic space. Remarkably, the existence of such a
collection of paths implies that the space is hyperbolic and the paths
are (reparametrized) quasi-geodesics with uniform constants. See
\cite{masur-schleimer}, which proves that arc complexes are
hyperbolic, and \cite{bowditch:criterion}.

Since the original proof of hyperbolicity of $\mathcal C(S)$ there
have been others, the simplest being \cite{unicorn}, not using
Teichm\"uller theory at all but constructing a family of paths in
$\mathcal C(S)$ directly using surgeries on curves. Perhaps
surprisingly, the more recent proofs also show that curve graphs are
\textit{uniformly} hyperbolic, i.e. $\delta$ can be taken
independently of the surface.

\subsection{The boundary of the curve complex}
If $X$ is a hyperbolic space which is not proper, its
boundary $\partial X$ may not be compact. For example, the boundary of
the wedge of countably many rays joined at the initial point is a
discrete countable set, and the boundary of a tree all of whose
vertices have countable valence is homeomorphic to the irrationals.

In
\cite{klarreich} E. Klarreich identified the boundary $\partial
\mathcal C$ of the curve complex as a proper quotient of a subspace of
Thurston's boundary of Teichm\"uller space $\mathcal T$. This
description serves as a model for boundaries of other hyperbolic
complexes.

\subsection{WPD, acylindrically hyperbolic groups, quasi-morphisms}

In the absence of properness of the action, one needs some kind of a
substitute. The property WPD (for ``weak proper discontinuity'') was
introduced in \cite{bf:wpd}.

\begin{definition}
  Suppose a group $G$ acts by isometries on a hyperbolic space
  $X$. A loxodromic element $g\in G$ is WPD if for every $x\in X$ and
  $C>0$ there is $N>0$ such that the set
  $$\{h\in G\mid d(x,h(x))\leq C, d(g^N(x),hg^N(x))\leq C\}$$
  is finite. The action of $G$ on $X$ is WPD if $G$ is not virtually
  cyclic and every loxodromic element is WPD.
\end{definition}

The WPD condition says that the collection of translates of an axis
(or an orbit) of a loxodromic element is discrete: any two translates
are either parallel or else they are in a bounded Hausdorff
neighborhood of each other only along a bounded length interval.

\begin{theorem}\cite{bf:wpd}\label{WPD}
  The action of $Mod(S)$ on $\mathcal C(S)$ is WPD. If a non-virtually
  cyclic group acts isometrically on a hyperbolic space with a
  WPD element then the space $\widetilde{QH}(G)$ of (reduced)
  quasi-morphisms on $G$ is infinite dimensional.
\end{theorem}

A \textit{quasi-morphism} is a function $f:G\to\R$ such that
$$\sup_{a,b\in G} |f(ab)-f(a)-f(b)|<\infty$$
Denote by $QH(G)$ the vector space of all quasi-morphisms on $G$ and
note the vector subspaces $Hom(G,\R)$ of homomorphisms $G\to\R$ and
$B(G)$ of bounded functions on $G$. Then the space $\widetilde{QH}(G)$
is defined as the quotient
$$\widetilde{QH}(G)=QH(G)/(Hom(G,\R)+B(G))$$
and it can also be identified with the kernel of the natural
homomorphism $H^2_b(G;\R)\to H^2(G;\R)$ from bounded cohomology of
$G$. For more on bounded cohomology see \cite{danny}. 
  
The basic method for showing $\widetilde{QH}(G)$ is
infinite-dimensional is due to Brooks \cite{brooks} in the case of
free groups. Fix a free group $F$ with a basis $a_1,a_2,\cdots$. Let
$w$ be any cyclically reduced word in the basis. Define
$f_w:F\to\Z\subset\R$ as $f_w(x)=C_w(x)-C_{w^{-1}}(x)$, where $C_w(x)$
is the number of occurrences of $w$ as a subword of $x$, written as a
reduced word. That $f_w$ is a quasi-morphism can be seen by
considering the tripod in the Cayley tree of $F$ spanned by $1,a$ and
$ab$ and marking all occurrences of $w^{\pm 1}$ along it. All such
occurrences that don't contain the central vertex will be counted
twice, with opposite signs, in the expression $f(ab)-f(a)-f(b)$, and of
course the number occurrences that do contain the central vertex is
uniformly bounded. With a bit more work one can show that for a
suitable choice of $w_i$'s the quasi-morphisms $f_{w_i}$ will yield
linearly independent elements of $\widetilde{QH}(F)$. The proof of the
second half of Theorem \ref{WPD} is a coarse version of this method,
where $w$ is replaced by a long segment along an axis of a WPD
element, and the discreteness of the set of translates guarantees that
the counting function is finite.

A quick application is the following statement, suggesting that
pseudo-Anosov elements of $Mod(S)$ are ``generic''.

\begin{corollary}
  Fix a finite generating set and the corresponding word metric on $Mod(S)$.
  For any $R>0$ there exists $M>0$ such that every ball of radius $M$
  contains a ball of radius $R$ that consist
  entirely of pseudo-Anosov mapping classes.
\end{corollary}

This follows quickly from the feature of the quasi-morphisms on
$Mod(S)$ constructed above that they are uniformly bounded on all
elements of $Mod(S)$ which are not pseudo-Anosov.

Bowditch noticed that the action of $Mod(S)$ on $\mathcal C(S)$
satisfies a property stronger than WPD.

\begin{definition}
  An isometric action of $G$ on a hyperbolic space $X$ is
  \textit{acylindrical} if for all $r>0$ there exist $R,N>0$ so that
  whenever $a,b\in X$ with $d(a,b)\geq R$, then there are at most $N$
  elements $h$ of $G$ such that $d(a,h(a))\leq r$ and $d(b,h(b))\leq
  r$.
\end{definition}

Thus acylindricity gives control in all directions, not only along
axes of loxodromic elements.

\begin{theorem}\cite{bowditch:acylindrical}
  The action of $Mod(S)$ on $\mathcal C(S)$ is acylindrical.
\end{theorem}

These results motivated Denis Osin to propose acylindrically
  hyperbolic groups as a generalization of hyperbolic groups. A group
is \textit{acylindrically hyperbolic} if it is not virtually cyclic and admits
an acylindrical action on a hyperbolic space with unbounded
orbits. This class contains many groups of interest (e.g. mapping
class groups and $Out(F_n)$) and many constructions on hyperbolic
groups carry over to this larger class, e.g. small cancellation
theory, or quasi-morphisms indicated above. See \cite{osin:ah,osin:icm}.

\subsection{Subsurface projections}

The main drawback of acylindrically hyperbolic groups is that in
general one doesn't have access to elements that don't act
loxodromically. In the case of mapping class groups this problem is resolved
through \textit{subsurface projections} of Masur and Minsky
\cite{MM1,MM2}.

Let $S$ be a surface as before and $X\subset S$ a connected
$\pi_1$-injective subsurface which is closed as a subset. Let $\alpha$
be a simple closed curve in $S$ which cannot be homotoped in the
complement of $X$ and which is in minimal position with respect to
$\partial X$. Then the intersection $\alpha\cap X$ consists of
finitely many disjoint arcs (or just $\alpha$ if $\alpha\subset
X$). For each such arc $J$ consider one or two curves obtained as
follows. If the endpoints of $J$ are contained in the same boundary
component $b$ of $X$, there are two ways of closing up $J$ to a closed
curve by adding an arc in $b$; take both of these curves. If the
endpoints of $J$ are on distinct boundary components $b,b'$ then form
a curve by taking two parallel copies of $J$ and connect them by
adding ``long'' arcs in $b$ and $b'$. It is not hard to see
that taking the collection of all these curves for all arcs $J$
produces a uniformly bounded set $\pi_X(\alpha)\subset\mathcal C(X)$
(we collapse all boundary components of $X$ to punctures). This
construction makes sense whenever $\mathcal C(X)$ is defined (so
notably a pair of pants is excluded). It also makes sense when $X$ is
an annulus, in which case the curve complex is formed by arcs joining
the boundary components, but we will not describe this case in
detail. If $\alpha$ is disjoint from $X$ then $\pi_X(\alpha)$ is not
defined and we set it to be empty.

Now fix a finite collection of curves $\vec{\alpha}=\{\alpha_1,\cdots,\alpha_n\}$ in
$S$ that ``fill'' the surface, i.e. every (essential) curve intersects
at least one of them. By the classical fact that the distance in the
curve complex is bounded by a function of the intersection number, if
$\pi_X(\alpha_i)$ and $\pi_X(\alpha_j)$ are both defined then their
union has uniformly bounded diameter (with the bound depending on the
intersection number between $\alpha_i$ and $\alpha_j$). We then
define $$\pi_X(\vec{\alpha})=\cup_i\pi_X(\alpha_i)$$
This is always a nonempty uniformly bounded subset of $\mathcal C(X)$.

The following is the fundamental result of Masur and Minsky,
expressing (coarsely) the word metric in $Mod(S)$ in terms of
subsurface projections. For $K>0$ and $x\geq 0$ define $\{\{x\}\}_K$
as 0 if $x<K$ and as $x$ if $x\geq K$.

\begin{theorem}[The distance formula]\cite{MM2}
  For all sufficiently large $K$ (depending on $\vec\alpha$) and for
  all $g,h\in Mod(S)$ we have
  $$d(g,h)\asymp \sum_X \{\{d_X(g(\vec\alpha),h(\vec\alpha))\}\}_K$$
\end{theorem}

The left hand side is the distance in the word metric. The summation
is over all (isotopy classes of) connected, $\pi_1$-injective
subsurfaces $X$ with $\mathcal C(X)\neq\emptyset$, and the displayed summand is
the diameter of the set $\pi_X(g(\vec\alpha))\cup
\pi_X(h(\vec\alpha))$. The symbol $\asymp$ means that there is a
linear function (depending on $K$ and the finite generating set of
$Mod(S)$) $f(x)=Ax+B$ such that the left-hand side is bounded by the
$f$-value of the right-hand side and vice versa. In particular, only
finitely many terms are $\geq K$.

The distance formula is a powerful tool in the study of large-scale
geometry of mapping class groups. It is used in an essential way in
the following remarkable theorem, establishing quasi-isometric
rigidity of mapping class groups. To state the theorem, let
$Mod^{\pm}(S)$ denote the \textit{extended} mapping class group,
i.e. allowing orientation-reversing homeomorphisms (this is an index 2
extension of $Mod(S)$). If $G$ is a finitely generated group with a
word metric, denote by $QI(G)$ the group of quasi-isometries $G\to G$
with the equivalence relation $f_1\sim f_2$ if $\sup_g
d(f_1(g),f_2(g))<\infty$. There is a natural homomorphism $G\to QI(G)$
sending $g$ to the left translation by $g$.

\begin{theorem}\cite{qi-rigidity,ursula2}
  Let $S$ be a surface of finite type. Except for a small number of
  sporadic surfaces, the natural homomorphism $Mod^{\pm}(S)\to QI(Mod^{\pm}(S))$
  is an isomorphism. In particular if $G$ is any group quasi-isometric
  to $Mod(S)$ there is a homomorphism $G\to Mod^{\pm}(S)$ with
  finite kernel and finite index image.
\end{theorem}

\section{Projection complexes}

It is tempting to view the distance formula as saying that the coarse map
$$Mod(S)\to\prod_X \mathcal C(X)$$ defined by $g\mapsto
\pi_X(g(\vec\alpha))$ is a quasi-isometric embedding, where we equip
the right-hand side with the $\ell_1$-metric. The trouble is that this
is not really a metric, and ``cutting off'' at $K$ in each coordinate
wouldn't satisfy the triangle inequality. Up to modifying each
coordinate a bounded amount, the image of this map was identified in
\cite{behrstock,qi-rigidity}. The main restriction on the image is the
following inequality.

\begin{theorem}[Behrstock inequality]\cite{behrstock}
  There is a $\theta\geq 0$ such that the following holds.
  Suppose $X,Y\subset S$ are two subsurfaces such that the boundary of
  each intersects the other. Then at least one of $d_X(\partial
  Y,\vec\alpha)$ and $d_Y(\partial
  X,\vec\alpha)$ is $\leq\theta$.
\end{theorem}

There is a simple proof of the Behrstock inequality, due to Chris
Leininger, see \cite{mangahas}. If we focus on the two coordinates
$\mathcal C(X)\times \mathcal C(Y)$, the inequality says that the
image is contained in a Hausdorff neighborhood of the ``wedge'' of
$\mathcal C(X)\times \{y\}\cup \{x\}\times \mathcal C(Y)$ where
$x=\pi_X(\partial Y)$ and $y=\pi_Y(\partial X)$. This suggests taking
wedges instead of products for the right-hand side in order to fix the
metrizability problem, and leads to the following construction that
can be axiomatized.

Let $\mathcal Y$ be a collection of metric spaces (technically we
allow the distance to be infinite, for example we might have
disconnected graphs with the the path metric). Suppose that for
distinct $X,Y\in\mathcal Y$ we are given a subset $\pi_X(Y)\subset
X$. If $Z\in\mathcal Y$, $Z\neq X$, we define
$$d_X(Y,Z)=\diam (\pi_X(Y)\cup\pi_X(Z))$$
We will assume that the following axioms hold for some fixed
$\theta\geq 0$:

\begin{enumerate}
\item [(P1)] $d_X(Y,Y)\leq\theta$,
\item [(P2)] if $d_X(Y,Z)>\theta$ then $d_Y(X,Z)\leq\theta$, and
  \item [(P3)] for $X\neq Z$ the set $$\{Y\in\mathcal Y\mid
    d_Y(X,Z)>\theta\}$$
    is finite.
\end{enumerate}

There are many natural situations where these axioms hold.

\begin{examples}
  \begin{enumerate}[(1)]
    \item Let $T$ be a simplicial tree and $\mathcal Y$ a collection
      of pairwise disjoint simplicial subtrees. The projection $\pi_X(Y)$
      is the point of $X$ nearest to $Y$. The axioms hold with
      $\theta=0$.

      \begin{figure}[h]
        \begin{center}
          \begin{tikzpicture}[scale=0.5,y=-1cm]
            \draw[black] (8.6,4.8) -- (8.6,9.1);
            \node at (9.2,4.2) {$\scriptstyle\pi_A(B)=\pi_A(C)$};
            \node at (16.3,4.2) {$\scriptstyle\pi_B(A)=\pi_B(C)$};
            \node at (8.6,9.7) {$\scriptstyle\pi_C(A)$};
            \node at (17.1,9.7) {$\scriptstyle\pi_C(B)$};
            \node at (6.1,4) {$A$};
            \draw[black] (17.1,4.8) -- (17.1,9.1);
            \node at (20.5,4.2) {$B$};
            \node at (19.5,10.3) {$C$};
            \draw [fill] (8.6,4.8) circle [radius=0.1];
            \draw [fill] (8.6,9.1) circle [radius=0.1];
            \draw [fill] (17.1,4.8) circle [radius=0.1];
            \draw [fill] (17.1,9.1) circle [radius=0.1];
\draw[very thick,black] (4.7,4.8) -- (9.4,4.8);
\draw[very thick,black] (7.2,2.7) -- (7.2,4.8);
\draw[very thick,black] (16.1,4.8) -- (20.8,4.8);
\draw[very thick,black] (18.2,2.7) -- (18.2,4.8);
\draw[very thick,black] (3.3,9.1) -- (21.6,9.1);
\draw[very thick,black] (5.7,9.1) -- (5.7,13.2);
\draw[very thick,black] (15.7,9.1) -- (15.7,13.2);

\end{tikzpicture}%
          \caption{The situation of Example (1), $d_C(A,B)>0$ while $d_A(B,C)=d_B(A,C)=0$.}
        \end{center}
        \end{figure}
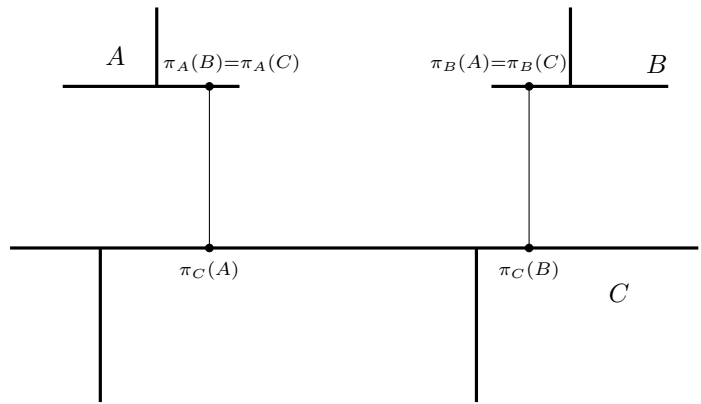
      
    \item Let $S$ be a closed hyperbolic surface and $\gamma$ an
      immersed closed geodesic which is not a multiple. In the
      universal cover $\tilde S=\mathbb H^2$ consider the set
      $\mathcal Y$ of all lifts of $\gamma$, and define projections as
      nearest point projections. A similar construction can be
      performed with a group acting on a hyperbolic space and a
      maximal virtually cyclic subgroup that contains a WPD element.
      \item Let $S$ be a complete hyperbolic surface of finite area
        and a cusp. In the universal cover $\tilde S=\mathbb H^2$
        consider the set $\mathcal Y$ of all lifts of a fixed horocyclic
        curve in the cusp (with either the intrinsic or the induced
        metric). Again the projection is the nearest point
        projection. A similar construction can be performed with
        relatively hyperbolic groups.
      \item Let $G$ be a group acting on a simplicial
        hyperbolic graph $X$ and let $H$ be the stabilizer of a vertex
        $v\in X$. Assume that $H$ acts simply transitively on the
        edges incident to $v$, and that the metric on the link
        $Lk(v,X)$ (which can be identified with $H$) induced by the
        path metric on $X\smallsetminus\{v\}$ is proper (finite radius
        balls contain finitely many points; here we allow distances to
        be infinite). Let $\mathcal Y$ be the collection of links of
        vertices in the orbit of $v$
        with this proper metric on each. If $u,w$ are two distinct
        vertices in the orbit of $v$ the projection of $Lk(u,X)$ to
        $Lk(w,X)$ is the set of points in $Lk(w,X)$ that belong to a
        geodesic between $u$ and $w$. If $(G,H)$ admit such an action,
        $H$ is said to be \textit{hyperbolically embedded} in $G$. See
        \cite{dgo}. For example, parabolic subgroups of hyperbolic
        groups, or maximal virtually cyclic subgroups containing a WPD
        element as in (2) are hyperbolically embedded, as can be seen
        by building the projection complex below.
     \item Let $S$ be an orientable surface of finite type and let
       $\mathcal Y$ be a collection of isotopy classes of
       $\pi_1$-injective subsurfaces where subsurface projections are
       defined, and assume that if $X,Y\in\mathcal Y$ and $X\neq Y$
       then $\partial X$ is not disjoint from $Y$ (up to
       isotopy). Define $\pi_Y(X)=\pi_Y(\partial X)$.
  \end{enumerate}
\end{examples}

The construction of a projection complex $\mathcal P(\mathcal Y)$ (and
the blow-up version $\mathcal C(\mathcal Y)$) is kind of a converse to
Example (2) above, where one tries to ``reconstruct'' the ambient
space from the projection data (though usually one gets a different
ambient space).

\begin{theorem}[\cite{bbf}, for a simpler construction see
    \cite{bbfs}]\label{pc} 
  Suppose the projection data $$(\mathcal Y,\pi_X(Y),\theta)$$ satisfy
  (P1-P3). There is a metric space $\mathcal C(\mathcal Y)$ containing
  metric spaces in $\mathcal Y$ as pairwise disjoint isometrically
  embedded subspaces and so that $\pi_X(Y)$ agrees, up to a bounded
  error, with the nearest point projection of $Y$ to $X$ within
  $\mathcal C(\mathcal Y)$. Moreover:
  \begin{itemize}
    \item If each $Y\in\mathcal Y$ is $\delta$-hyperbolic for some
      $\delta\geq 0$ then $\mathcal C(\mathcal Y)$ is 
      hyperbolic.
      \item If each $Y\in \mathcal Y$ is quasi-isometric to a tree (a
        ``quasi-tree'') 
        with fixed QI constants, then $\mathcal C(\mathcal Y)$ is also
        a quasi-tree.
        \item If the collection $\mathcal Y$ consists of finitely
          many isometry types of metric spaces and they all have
          asymptotic dimension $\leq n$ then $\asdim \mathcal
          C(\mathcal Y)\leq n+1$.
          \item The space $\mathcal P(\mathcal Y)$ obtained from
            $\mathcal C(\mathcal Y)$ by collapsing all embedded copies
            of spaces in $\mathcal Y$ is a quasi-tree.
            \item If a group $G$ acts by isometries on $\sqcup_{Y\in
              \mathcal Y}Y$ preserving the projections
              (i.e. $g(\pi_X(Z)=\pi_{g(X)}(g(Z))$ for all $g\in G$)
              then $G$ acts by isometries on $\mathcal C(\mathcal Y)$
              extending the action on $\sqcup_{Y\in
              \mathcal Y}Y$, and it also acts isometrically on
              $\mathcal P(\mathcal Y)$. 
  \end{itemize}
\end{theorem}

We briefly outline the construction. As indicated above, the idea is
to start with the disjoint union of all $Y\in\mathcal Y$ and then for
certain pairs $(X,Z)$ add edges joining points in $\pi_X(Z)$ to points
in $\pi_Z(X)$.

Step 1 is to promote (P2) to a stronger property (P2++):

\begin{enumerate}
\item [(P2++)] If $d_Y(X,Z)>\theta$ then $\pi_Y(X)=\pi_Y(Z)$.
\end{enumerate}

This can be done by modifying the projection $\pi_X(Y)$ by a bounded
amount and replacing $\theta$ by a larger constant. This modification
preserves group equivariance.

In step 2, assuming (P1,P2++,P3), one chooses a constant $K\geq
2\theta$ and posits that $X$ and $Z$ are connected by edges as above
provided $d_Y(X,Z)\leq K$ for all $Y\neq X,Z$. The key property that
makes the proof of Theorem \ref{pc} possible is that the set
$$\{X\}\cup \{Y\mid d_Y(X,Z)>K\}\cup \{Z\}$$ is finite (by (P3)) and
is naturally linearly ordered giving a path from $X$ to $Z$, called a
\textit{standard path}, in $\mathcal P_K(\mathcal Y)$. These standard
paths are quasi-geodesics and behave very nicely. The construction
depends on the choice of the constant $K$: when $K$ is enlarged, there
will be more edges attached.

We mention a few applications of this construction to mapping class
groups.

\begin{theorem}\cite{bbf}\label{mcg-asdim}
  $\asdim(Mod(S))<\infty$.
\end{theorem}

The basic idea is to replace the infinite product of curve complexes
by a smaller space. The collection of all subsurfaces $\mathcal Y$
does not satisfy the assumptions of Example (5) above since
subsurfaces can be disjoint or nested. However, one shows that there
is a way to write $\mathcal Y$ equivariantly as a finite disjoint union
$\sqcup\mathcal Y_i$ so that each collection $\mathcal Y_i$ satisfies
Example (5). Thus one gets the spaces $\mathcal C(\mathcal
Y_i)$. These are all hyperbolic, and crucially, have finite asymptotic
dimension by Theorem \ref{pc} and the theorem of Bell-Fujiwara
\cite{bell-fujiwara} that curve complexes have finite asymptotic
dimension. Then we have a quasi-isometric embedding $$Mod(S)\to\prod_i
\mathcal C(\mathcal Y_i)$$
which finishes the proof since passing to finite products and
subspaces preserves finiteness of asymptotic dimension.

There is quite a bit of inefficiency when we take the product of the
blown-up projection complexes over the families $\mathcal Y_i$. There
is a more involved system of axioms that keeps track of pairs of
surfaces that are disjoint or nested leading to the notion of a
\textit {hierarchically hyperbolic group}, due to J. Behrstock,
M. Hagen and A. Sisto. For example, in \cite{hhs:asdim} they derive a
bound on $\asdim(Mod(S))$, using \cite{bb:asdim}, which is quadratic
in the complexity of the surface. There are other applications of this
theory, for example in \cite{hhs:qf} they show how to understand
quasi-flats in mapping class groups and how to approximate a ``hull''
of a finite set by a $CAT(0)$ cube complex.

\begin{theorem}\cite{bbf:scl}
  There is a classification, in terms of the Nielsen-Thurston normal
  form, of those elements $g$ of $Mod(S)$ that have stable commutator
  length $scl(g)=0$.
\end{theorem}

Recall that for $g\in [G,G]$ $cl(g)$ is the smallest $k$ such that $g$
can be written as a product of $k$ commutators, and
$scl(g)=\lim_n\frac{cl(g^n)}n$. By Bavard duality (see \cite{danny}),
$scl(g)>0$ is equivalent to having a quasi-morphism $G\to\R$ which is
unbounded on the powers of $g$. Projection complexes are used to
construct actions of finite index subgroups of $Mod(S)$ on hyperbolic
spaces with a power of a given element acting loxodromically, and then
the Brooks method can be used to construct such quasi-morphisms. It's
worth stating this fact:

\begin{theorem}\cite{bbf:scl}\label{BBF}
  Let $S$ be a finite type surface. There is a torsion-free finite
  index subgroup $G<Mod(S)$ such that for every element $g\in G$
  of infinite order there an
  action of $G$ on a hyperbolic space such that $g$
  is loxodromic.
\end{theorem}

For example, this applies to (powers of) Dehn twists. By contrast, a
theorem of Bridson \cite{bridson:twist} says that whenever $Mod(S)$
(with $S$ of genus $\geq 3$) acts on a $CAT(0)$ space, Dehn twists
have translation length 0. 

Projection complexes are useful more generally for constructing
quasi-cocycles on groups $G$ with coefficients in orthogonal
representations on strictly convex Banach spaces (such as $l^p(G)$ for
$1<p<\infty$). See \cite{bbf:banach}.

\begin{theorem}[Balasubramanya \cite{sahana}]
  If a group $G$ acts on a hyperbolic space with a WPD element, then
  it admits a cobounded acylindrical action on a quasi-tree.
\end{theorem}

Another proof of Balasubramanya's theorem is given in \cite{bbfs}. The
quasi-tree is the projection complex applied to Example (2) and
acylindricity is proved using the geometry of standard paths.

F. Dahmani, V. Guirardel and D. Osin solved a long standing open
problem when they proved the following.

\begin{theorem}\cite{dgo}\label{pA}
  Let $\phi\in Mod(S)$ be a pseudo-Anosov mapping class. Then for a
  suitable power $\phi^n$ with $n>0$ the subgroup normally generated
  by $\phi^n$ is free.
\end{theorem}

They derive this theorem using the method of \textit{rotating
families}.

\begin{theorem}\cite{dgo}\label{rotation}
  For every $\delta\geq 0$ there is $R>0$ such that the following holds.
  Let $X$ be a $\delta$-hyperbolic space and $G$ a group of isometries
  of $X$. Let $C\subset X$ be a $G$-invariant set which is
  $R$-separated (meaning that $d(c,c')>R$ if $c,c'\in C$ are
  distinct). Suppose for every $c\in C$ we are given a subgroup $G_c$
  of the stabilizer $Stab_G(c)$ such that
  \begin{enumerate}[(i)]
  \item $G_{g(c)}=gG_cg^{-1}$ for $c\in C$ and $g\in G$, and
    \item if $g\in G_c\smallsetminus\{1\}$, $c'\in C$ and $c'\neq c$ then every
      geodesic from $c'$ to $g(c')$ passes through $c$.
  \end{enumerate}
  Then the subgroup of $G$ generated by $\cup_{c\in C} G_c$ is the free
  product of a subcollection of the family $\{G_c\}_{c\in C}$.
\end{theorem}

To prove Theorem \ref{pA} they apply this theorem to the space
obtained from the curve complex $\mathcal C(S)$ by equivariantly
coning off an orbit of the elementary closure $EC(\phi)$. Pretending
that this orbit is in an isometrically embedded line, one would attach
the universal cover of a disk of large radius in $\mathbb H^2$ punctured at the
center, and then completed to add the cone point back in. The set of
these cone points is the set $C$ from the theorem, and $G_c$ is the
cyclic group generated by (a conjugate of) a suitable power $\phi^n$.

More recently, M. Clay, J. Mangahas and D. Margalit proved
a version of Theorem \ref{rotation} that applies to projection
complexes. Rotating families are replaced by \textit{spinning families}.

\begin{theorem}\cite{cmm}\label{CMM}
For every $\theta$ and $K$ there is $L$ so that the following
holds. Suppose a group $G$ acts on the projection data and on the
associated projection complex $\mathcal P=\mathcal P_K(\mathcal
Y)$. Suppose for every vertex $v\in\mathcal P$ we are given a subgroup
$G_v$ of the stabilizer $Stab_G(v)$ such that
\begin{enumerate}[(i)]
  \item $G_{g(v)}=gG_vg^{-1}$ for any vertex $v$ and $g\in G$, and
\item if $v,v'$ are distinct vertices and $g\in
  G_v\smallsetminus\{1\}$ then $d_v(v',g(v'))>L$.
\end{enumerate}
Then the subgroup of $G$ generated by $\cup_{v} G_v$ is the free
  product of a subcollection of the family $\{G_v\}_{v\in \mathcal
    P^{(0)}}$.
\end{theorem}

They derive Theorem \ref{pA} directly from Theorem \ref{CMM} using the
projection complex as in Example (2). They also prove several
statements about normal closures of powers of other kinds of elements,
or collections of elements, in $Mod(S)$. One extreme behavior is that
the normal closure is free, another that it is the whole $Mod(S)$, but
surprisingly there are examples when the normal closure turns out to
be a certain kind of (infinitely generated) right angled Artin groups.

In \cite{manymany} the two theorems above are revisited
and in particular the paper shows how to derive Theorem \ref{rotation}
from Theorem \ref{CMM}.

Here are two more applications of projection complexes to mapping
class groups, though we will not comment on
the proofs.

\begin{theorem}\cite{bartels-b}
  Mapping class groups satisfy the Farrell-Jones conjecture.
\end{theorem}

\begin{theorem}\cite{combing}
  Mapping class groups are semi-hyperbolic.
\end{theorem}

This means that one can equivariantly choose uniform quasi-geodesics
connecting any pair of points in $Mod(S)$ so that they fellow-travel,
i.e. if the endpoints are at distance $\leq 1$ then each is in the
other's uniform Hausdorff neighborhood.

\section{$Out(F_n)$}

Let $F_n$ be the free group of rank $n\geq 2$, $Aut(F_n)$ its
automorphism group, and $Out(F_n)=Aut(F_n)/F_n$ the \textit{outer
  automorphism group} of $F_n$, obtained by quotienting out the inner
automorphisms. This group has been studied for over a century, see
Nielsen's paper \cite{nielsen} where he proves that $Out(F_n)$ is
generated by $n+1$ involutions. A big impediment in the study of
$Out(F_n)$, and free groups in general, was the tendency to think of
elements of free groups as words in a basis. A much more flexible
approach is to think of a free group as the fundamental group of a
graph, which is not necessarily a rose $R_n$ (a wedge of $n$ circles). For
example, the proof that subgroups of free groups are free is
essentially trivial using covering spaces and general graphs, while
the more algebraic proof is much less transparent. In
\cite{stallings:fold} J. Stallings introduced the operation of
\textit{folding} of graphs and used it to show that many standard
algorithmic problems about free groups have easy solutions.

\subsection{Outer space}
Given this philosophy, the definition of Culler-Vogtmann's Outer space
$CV_n$ should seem very natural. Fix the rose $R_n$. A point in $CV_n$
is represented by a homotopy equivalence $h:R_n\to\Gamma$, called
\textit{marking}, where
$\Gamma$ is a finite graph with all vertices of valence $>2$ equipped
with a \textit{metric} of volume 1, i.e. an assignment of positive
numbers to its edges that add to 1. Two such markings $h:R_n\to
\Gamma$ and $h':R_n\to\Gamma'$ represent the same point in $CV_n$ if
there is an isometry $\phi:\Gamma\to\Gamma'$ such that $\phi h$ is
homotopic to $h'$. Formally, the definition is analogous to the
definition of Teichm\"uller space, where metric graphs are replaced by
hyperbolic surfaces. There are many useful analogies between mapping
class groups and $Out(F_n)$, perhaps stemming from the classical
theorem of Dehn-Nielsen-Baer (see \cite{farb-margalit}) that when $G$
is the fundamental group of a closed orientable surface $S$ then
$Out(G)\cong Mod^{\pm}(S)$. While Teichm\"uller space is diffeomorphic
to Euclidean space, Outer space is a contractible polyhedron and the
study of $Out(F_n)$ is decidedly more combinatorial compared to the
study of mapping class groups. The group $Out(F_n)$ acts naturally on
$CV_n$ by changing the marking. The action is proper.
For more on Outer space and the
consequences to the structure of $Out(F_n)$ see the original paper
\cite{CV} as well as the excellent survey \cite{karen1}, and also
\cite{mb:icm}.

\subsection{The boundary of Outer space}
By taking universal covers, another way to think about a point
$h:R_n\to\Gamma$ in $CV_n$ is as a free action of $F_n$ on a
simplicial metric tree. The construction in section \ref{2.10} then
yields a compactification of $CV_n$ with the points in the ideal
boundary $\partial CV_n$ represented by actions of $F_n$ on $\R$-trees
(which are either non-simplicial or non-free). This construction was
carried out in \cite{culler-morgan}. Exactly which trees arise in
$\partial CV_n$ was identified in \cite{outerlimits,camille}.

\subsection{Lipschitz metric and train-track maps}
There is a natural notion of a \textit{Lipschitz distance} between two
points $h_i:R_n\to \Gamma_i$, $i=1,2$. It is defined
by $$d(\Gamma_1,\Gamma_2)=\log\lambda$$ where $\lambda\geq 1$ is the
smallest possible Lipschitz constant of all maps
$f:\Gamma_1\to\Gamma_2$ that commute with markings, i.e. $h_2 f$ is
homotopic to $h_1$ (and $\Gamma_i$ are viewed as geodesic metric
spaces). This ``metric'' is not symmetric, but satisfies the triangle
inequality $d(\Gamma_1,\Gamma_3)\leq
d(\Gamma_1,\Gamma_2)+d(\Gamma_2,\Gamma_3)$, and $d(\Gamma,\Gamma')\geq
0$ with equality only for $\Gamma=\Gamma'$. This metric has
interesting properties and displays a mixture of behaviors of the
well-studied metrics on Teichm\"uller space (Teichm\"uller,
Weil-Petersson and Thurston metrics). It can be used in the $Out(F_n)$
setting in a way similar
to the Bers' proof of the Nielsen-Thurston classification of mapping
classes (see \cite{farb-margalit}) to give a proof of the following
train-track theorem. See \cite{pcmi}. 

\begin{theorem}\cite{BH}
  Every irreducible automorphism $\phi\in Out(F_n)$ can be represented
  by a train-track map $f:\Gamma\to\Gamma$ for some $\Gamma\in CV_n$.
\end{theorem}

A marking gives an identification between $\pi_1(\Gamma)$ and $F_n$
and $f:\Gamma\to\Gamma$ ``represents'' $\phi$ if the induced
endomorphism on $\pi_1(\Gamma)$ is $\phi$. We say that $\phi$ is
\textit{irreducible} if it cannot be represented by some
$f:\Gamma\to\Gamma$ that leaves a proper subgraph with nontrivial
$\pi_1$ invariant. The map $f$ is a \textit{train-track map} if all
positive powers of $f$ are locally injective on all edges of
$\Gamma$. It is easy to control the growth of lengths of loops under
iteration by train-track maps, which makes them important in the study
of the dynamics of an automorphism. More generally, when $\phi$ is not
irreducible, there are \textit{relative} train-track representatives.

The Lipschitz metric admits geodesic paths, called \textit{folding
  paths}, which are induced, in the spirit of Stallings, by
identifying segments of the same length and issuing from the same
vertex. For more on this see \cite{pcmi}.   

\subsection{Hyperbolic complexes}

By analogy with the arc and curve complexes, there are several
complexes where $Out(F_n)$ acts.

\subsubsection{The free splitting complex $FS_n$}
This one is analogous to the arc complex.
A $k$-simplex is a $(k+1)$-edge free splitting of $F_n$, i.e. it is a
minimal action of $F_n$ on a simplicial tree with vertices of valence
$>2$, with trivial edge
stabilizers and with $(k+1)$ orbits of edges. Passing to a face is
induced by equivariantly collapsing an orbit of edges. Outer space
$CV_n$ is naturally a subset of $FS_n$, which can be viewed as a
``simplicial completion'' of $CV_n$.

\subsubsection{The cyclic splitting complex $FZ_n$}
This is defined the same way, except that the edge stabilizers can be
cyclic subgroups. It is analogous to the curve complex.

\subsubsection{The free factor complex $FF_n$}
This one is different from $FZ_n$ but can also be viewed as an analog
of the curve complex.
A vertex of $FF_n$ is a \textit{proper free factor} $A<F_n$, i.e. a
subgroup such that $F_n=A*B$ for some $A\neq 1\neq B$, defined up to
conjugation. A $k$-simplex is a $k$-tuple of distinct conjugacy
classes of proper free
factors that are nested after suitable conjugation.

There are natural coarse equivariant maps
$$CV_n\to FS_n\to FZ_n\to FF_n$$
For example, $FS_n\to FF_n$ sends a free splitting to a nontrivial
vertex group (or if they are all trivial to a free factor represented
by a subgraph of the quotient graph).

Now, it turns out that all three of these complexes are hyperbolic,
and there are several others that this survey is not mentioning. The
first hyperbolic $Out(F_n)$-complex was constructed in
\cite{bf:hyperbolic}, though it's not canonical. The hyperbolicity of
$FF_n$ was established in \cite{bf:FF} along the lines of the
Masur-Minsky's argument for the curve complex, by projecting folding
paths from $CV_n$ to $FF_n$. A novel argument by Handel-Mosher
\cite{handel-mosher} established hyperbolicity of $FS_n$, by
considering folding paths directly in $FS_n$. Kapovich-Rafi
\cite{ilya-kasra} found a general criterion that a
Lipschitz map
$X\to Y$ has to satisfy in order for the hyperbolicity of $X$ to imply
the hyperbolicity of $Y$. Essentially, Lipschitz images of thin
triangles are thin triangles. The maps $FS_n\to FZ_n\to FF_n$ satisfy
the Kapovich-Rafi criterion, so the hyperbolicity of $FS_n$ implies
the hyperbolicity of the other two. Loxodromic elements in $FF_n$ are
precisely the \textit{fully irreducible automorphisms} (the ones whose
positive powers are irreducible) and they are all WPD (in $FS_n$ there
are more loxodromic elements and they are not all WPD). Thus the space
of quasi-morphisms on $Out(F_n)$ is infinite-dimensional and
$Out(F_n)$ is acylindrically hyperbolic. Handel and Mosher
\cite{hm:qm1,hm:qm2} extended this and proved the $H^2_b$-alternative:
any subgroup of $Out(F_n)$ which is not virtually abelian has an
infinite-dimensional space of quasi-morphisms. This recovers the theorem
of Bridson and Wade \cite{bridson-wade} that no higher rank lattice
embeds as a subgroup of $Out(F_n)$. The proof is much more involved
than the $H^2_b$-alternative for mapping class groups \cite{bf:wpd}. 

The boundary of $FF_n$ was identified with a proper quotient of a
subspace of $\partial CV_n$ in \cite{boundaryFF} and in
\cite{ursula:boundary}.

\subsection{Subfactor projections}
By analogy with the Masur-Minsky subsurface projections, there are
\textit{subfactor projections}, see \cite{bf:proj, sam}.  Let $A,B$ be
two proper free factors in $F_n$. Our goal is to define $\pi_A(B)\in
FS(A)$, the projection of $B$ to the free splitting complex of
$A$. Choose $\Gamma\in CV_n$ so that $B$ is represented by a subgraph
$\Gamma_B$ of $\Gamma$. Then represent $A$ by an immersion
$\Gamma_A\to \Gamma$. Thus $\Gamma_A$ determines a simplex in Outer
space for $A$, and can be projected to $FS(A)$ (or $FF(A)$). It takes
some work to show that coarsely this projection does not depend on the
choice of $\Gamma$, at least when $A$ and $B$ are sufficiently far
apart in $FF_n$. Moreover, the set $\mathcal Y$ of all free factors
can be equivariantly and finitely partitioned into $\sqcup\mathcal
Y_i$ so that projection is defined within each $\mathcal Y_i$, and
this projection satisfies the projection axioms. One then gets a map
$$Out(F_n)\to \prod_i \mathcal C(\mathcal Y_i)$$
in the same way as for mapping class groups (see the discussion after
Theorem \ref{mcg-asdim}). However, here the map is \textit{not} a
quasi-isometric embedding. The main issue is that there is no analog
of annulus projections: when $A$ has rank 1, the corresponding complex
$FS(A)$ is a single point. For example, the orbits on the right-hand
side under the powers of any polynomially growing automorphism are
bounded. For more on this, see \cite{karen2}.

\subsection{Questions}
The following is the key question, if one hopes to understand
$Out(F_n)$ using hyperbolic methods. The other questions reiterate the
state of affairs that the large scale geometry of $Out(F_n)$ is
lagging behind the one of mapping class groups.

\begin{enumerate}
  \item [(1)] Given $\phi\in Out(F_n)$ of infinite order, is there a
    finite index subgroup $G<Out(F_n)$ and an isometric action of $G$
    on a hyperbolic space so that a positive power of $\phi$ that
    belongs to $G$ acts loxodromically?
\end{enumerate}

This is true for mapping class groups (see Theorem \ref{BBF}), and it is
also true for automorphisms $\phi$ that grow exponentially.

\begin{enumerate}
  \item [(2)] Do any of hyperbolic $Out(F_n)$-complexes admit
    \textit{tight (quasi-)geodesics}?
\end{enumerate}

These were defined for curve complexes by Masur and Minsky, and a very
strong finiteness property was established by Bowditch
\cite{bowditch:acylindrical}. Thus the question is asking for an
equivariant collection of uniform quasi-geodesics so that any two are
connected by at least one, but only finitely many of these.

Bowditch used his strong finiteness of tight geodesics to show that
translation lengths in the curve complex are rational, and
Bell-Fujiwara \cite{bell-fujiwara} used it to show that curve
complexes have finite asymptotic dimension.

\begin{enumerate}
  \item [(3)] Do the hyperbolic $Out(F_n)$-complexes $FS_n,FZ_n,FF_n$
    have finite asymptotic dimension? Are the translation lengths
    always rational? Does $Out(F_n)$ have finite
    asymptotic dimension?
\end{enumerate}

We remark that the Novikov conjecture is known for $Out(F_n)$,
\cite{bgh}.

The following seems out of reach with the present methods, although
\cite{pacman} is a promising start.

\begin{enumerate}
\item [(4)] Does $Out(F_n)$ satisfy the Farrell-Jones conjecture?
  \item [(5)] Does the local and global connectivity of $\partial
    FF_n$ go to infinity as $n\to\infty$? 
\end{enumerate}

By the work of Gabai \cite{gabai} the answer is yes for the boundary
of the curve complex. Each $\partial FF_n$ is finite dimensional
\cite{rick} and \cite{bch} is a start. Of course, the same question
can be asked about the boundaries of $FZ_n$ and $FS_n$.

\bibliographystyle{plain}
\bibliography{ref}

\end{document}